\documentclass[11pt, reqno, psamsfonts]{amsart}
\pdfoutput=1
\usepackage{amssymb}
\usepackage{amsthm}
\usepackage{amsmath}
\usepackage{latexsym}
\usepackage[T1]{fontenc}
\usepackage[utf8]{inputenc}
\usepackage[russian, french, english]{babel}
\usepackage{graphicx}
\usepackage{wrapfig}
\usepackage[justification=centering, labelfont=bf]{caption}
\usepackage{subcaption}
\usepackage{mathtools}
\usepackage[hidelinks]{hyperref}
\usepackage{enumitem}
\usepackage{mathrsfs}
\usepackage{thmtools}
\usepackage{stmaryrd}
\usepackage{framed}
\usepackage{tikz}
\usepackage{mathabx}
\usepackage{float}
\usepackage{bbm}
\usepackage{titlesec}

\usepackage{lmodern}

\usepackage[shortcuts]{extdash}

\usepackage{amsaddr}

\usetikzlibrary{shapes,snakes}
\usetikzlibrary{arrows.meta}

\usepackage[spacing=true,kerning=true,babel=true,tracking=true]{microtype}

\usepackage[backend=biber, style=alphabetic, sorting=nyt, maxnames=100]{biblatex}
\title{DP-Colorings of Hypergraphs}
\date{}
\author{Anton~Bernshteyn}
\address[Anton Bernshteyn]{Department of Mathematics, University of Illinois at Urbana--Champaign, IL, USA}
\email{bernsht2@illinois.edu}

\author{Alexandr Kostochka}
\address[Alexandr Kostochka]{Department of Mathematics, University of Illinois at Urbana--Champaign, IL, USA and Sobolev Institute of Mathematics, Novosibirsk, Russia}
\email{kostochk@math.uiuc.edu}
\thanks{Research of the first author is supported in part by the Waldemar J., Barbara G., and Juliette Alexandra Trjitzinsky Fellowship. Research of the second author is supported in part by NSF grant
DMS-1600592 and grants 18-01-00353A and 16-01-00499 of the Russian Foundation for Basic Research.}

\usepackage[left=1in,right=1in,top=1in,bottom=1in,bindingoffset=0cm]{geometry}

\newtheoremstyle{bfnote}%
{}{}%
{\slshape}{}%
{\bfseries}{\bfseries.}%
{ }%
{\thmname{#1}\thmnumber{ #2}\thmnote{ \ep{\normalfont{}#3}}}

\newtheoremstyle{defsc}%
{}{}%
{}{}%
{\scshape}{.}%
{ }%
{\thmname{#1}\thmnumber{ #2}\thmnote{ (#3)}}

\newtheoremstyle{claim}%
{}{}%
{\slshape}{}%
{\itshape}{.}%
{ }%
{\thmname{#1}\thmnumber{ #2}\thmnote{ \ep{\normalfont{}#3}}}

\theoremstyle{bfnote}

\numberwithin{equation}{section}

\newtheorem{theo}[equation]{Theorem}

\newtheorem{lemma}[equation]{Lemma}
\newtheorem{corl}[equation]{Corollary}

\newtheorem*{claim*}{Claim}
\newtheorem*{corl*}{Corollary}

\theoremstyle{definition}
\newtheorem{defn}[theo]{Definition}

\newtheorem{prob}[theo]{Problem}
\newtheorem{remk}[equation]{Remark}

\newtheorem{exmp}[theo]{Example}

\theoremstyle{remark}
\newtheorem*{remk*}{Remark}

\newcommand*{\myproofname}{Proof}
\newenvironment{claimproof}[1][\myproofname]{\begin{proof}[#1]}{\end{proof}}

\newcommand{\0}{\varnothing}
\newcommand{\set}[1]{\{#1\}}
\newcommand{\dom}{\mathrm{dom}}

\newcommand{\N}{\mathbb{N}}
\newcommand{\Z}{\mathbb{Z}}
\newcommand{\R}{\mathbb{R}}
\let\originalH\H
\usepackage{xspace}
\newcommand{\Erdos}{Erd\originalH os\xspace}
\renewcommand{\H}{\mathcal{H}}
\newcommand{\F}{\mathcal{F}}
\newcommand{\Q}{\mathbb{Q}}

\renewcommand{\epsilon}{\varepsilon}
\renewcommand{\phi}{\varphi}
\renewcommand{\theta}{\vartheta}
\renewcommand{\tilde}{\widetilde}
\renewcommand{\leq}{\leqslant}
\renewcommand{\geq}{\geqslant}

\newcommand{\defeq}{\coloneqq}

\newcommand{\fun}[2]{[#1 \rightharpoonup #2]}
\newcommand{\emphd}[1]{{\emph{#1}}}
\renewcommand{\P}{\mathbb{P}}

\newenvironment{coolproof}[1][Proof]{\begin{proof}[\textsc{#1}]}{\end{proof}}

\titleformat{\section}[block]{\scshape\filcenter}{\thesection.}{1ex}{}
\titleformat{\subsection}[block]{\bfseries\filcenter}{\thesubsection.}{1ex}{}
\titleformat{\subsubsection}[runin]{\bfseries}{\thesubsubsection.}{1ex}{}[.]

\titlespacing*{\section}{0pt}{*3}{*1}
\titlespacing*{\subsection}{0pt}{*3}{*1}

\makeatletter
\newcommand{\neutralize}[1]{\expandafter\let\csname c@#1\endcsname\count@}
\makeatother

\newenvironment{theocopy}[1]
{\renewcommand{\thetheo}{\ref{#1}}
	\neutralize{theo}\phantomsection
	\begin{theo}}
	{\end{theo}}

\newcommand{\bemph}[1]{{\upshape#1}} % define how emphasised brackets should look
\newcommand{\ep}[1]{\bemph{(}#1\bemph{)}} % parentheses

\renewcommand{\mkbibnamefirst}[1]{\textsc{#1}}
\renewcommand{\mkbibnamelast}[1]{\textsc{#1}}
\renewcommand{\mkbibnameprefix}[1]{\textsc{#1}}
\renewcommand{\mkbibnameaffix}[1]{\textsc{#1}}
\renewcommand{\finentrypunct}{}

\renewbibmacro{in:}{}

\renewbibmacro*{volume+number+eid}{%
	\printfield{volume}%
	\setunit*{\addnbspace}% NEW (optional); there's also \addnbthinspace
	\printfield{number}%
	\setunit{\addcomma\space}%
	\printfield{eid}}
\DeclareFieldFormat[article]{volume}{\textbf{#1}\space}
\DeclareFieldFormat[article]{number}{\mkbibparens{#1}}

\DeclareFieldFormat{journaltitle}{#1,}
\DeclareFieldFormat[thesis]{title}{\mkbibemph{#1}\addperiod}
\DeclareFieldFormat[article, unpublished, thesis]{title}{\mkbibemph{#1},}
\DeclareFieldFormat[book]{title}{\mkbibemph{#1}\addperiod}
\DeclareFieldFormat[unpublished]{howpublished}{#1, }

\DeclareFieldFormat{pages}{#1}

\DeclareFieldFormat[article]{series}{Ser.~#1\addcomma}

\renewcommand*{\bibfont}{\small}

\begin{document}
	\maketitle
	
	\begin{abstract}
		Classical problems in hypergraph coloring theory are to estimate the minimum number of edges, $m_2(r)$ (respectively, $m^\ast_2(r)$), in a non-$2$-colorable $r$-uniform (respectively, $r$-uniform and simple) hypergraph. The best currently known bounds are
		\[
			c \cdot \sqrt{r/\log r} \cdot 2^r \,\leq\, m_2(r) \,\leq\, C \cdot r^2 \cdot 2^r \qquad \text{and} \qquad
			c' \cdot r^{-\epsilon} \cdot 4^r \,\leq\, m_2^\ast(r) \,\leq\, C' \cdot r^4 \cdot 4^r,
		\]
	for any fixed $\epsilon > 0$ and some $c$, $c'$, $C$, $C' > 0$  \ep{where $c'$ may depend on $\epsilon$}. In this paper we consider the same problems in the context of \emph{DP-coloring} \ep{also known as \emph{correspondence coloring}}, which is a generalization of list coloring introduced by Dvo\v r\'ak and Postle %\cite{DP}
	and related to \emph{local conflict coloring} studied independently by Fraigniaud, Heinrich, and Kosowski. %\cite{FHK}.
	Let $\tilde{m}_2(r)$ (respectively, $\tilde{m}^\ast_2(r)$) denote the minimum number of edges in a non-$2$-DP-colorable $r$-uniform (respectively, $r$-uniform and simple) hypergraph. By definition, $\tilde{m}_2(r) \leq m_2(r)$ and $\tilde{m}^\ast_2(r)\leq m^\ast_2(r)$.
	
	While the proof of the bound $m^\ast_2(r) = \Omega( r^{-3} 4^r)$ due to \Erdos and Lov\'asz also works for $\tilde{m}^\ast_2(r)$,	we show that the trivial lower bound $\tilde{m}_2(r) \geq 2^{r-1}$ is asymptotically tight, i.e., $\tilde{m}_2(r) \leq (1 + o(1))2^{r-1}$. %, in contrast to the bound $m_2(r) > a \cdot \sqrt{r/\log r} \cdot 2^r$.
	On the other hand, when $r \geq 2$ is even, we prove that the lower bound $\tilde{m}_2(r) \geq 2^{r-1}$ is not sharp, i.e., $\tilde{m}_2(r) \geq 2^{r-1}+1$. Whether this result holds for any odd values of $r$ remains an open problem. Nevertheless, we conjecture that the difference $\tilde{m}_2(r) - 2^{r-1}$ can be arbitrarily large.

 
 \end{abstract}
	
	\section{Introduction}
	
	\subsection{Hypergraph coloring}
	
	We use $\N$ to denote the set of all nonnegative integers and write $\N^+ \defeq \N \setminus\set{0}$. For convenience, each integer $k \in \N$ is identified with the $k$-element set $\set{i \in \N \,:\, i < k}$. A \emphd{hypergraph} $\H$ on a set $X$ of \emphd{vertices} is a collection of subsets of $X$, called the \emphd{edges} of $\H$. A hypergraph $\H$ is \emphd{$r$-uniform} for $r \in \N^+$ if $|e|=r$ for all $e \in \H$. An \emphd{independent set} in $\H$ is a subset $I \subseteq X$ such that $e \not\subseteq I$ for all $e \in\H$. For $k \in \N$, a \emphd{proper $k$-coloring} of $\H$ is a mapping $f \colon X \to k$ such that for each $0 \leq i < k$, the preimage $f^{-1}(i) \subseteq X$ is an independent set in $\H$.
	
	\begin{defn}[{\Erdos--Hajnal \cite{ErdHaj}}]\label{defn:ErdHaj}
		For $r$, $k \in \N^+$, let $m_k(r)$ denote the minimum $|\H|$ over all $r$-uniform hypergraphs $\H$ that do not admit a proper $k$-coloring.
	\end{defn}
	
	\noindent We will be particularly interested in the case $k = 2$. The known values of $m_2(r)$ are:
	\[
		m_2(1) = 1, \qquad m_2(2) = 3, \qquad m_2(3) = 7, \qquad \text{and} \qquad m_2(4) = 23.
	\]
	The last of these was found by \"Osterg\aa rd~\cite{Ost} through an exhaustive computer search.
	
	A simple first-moment argument due to \Erdos \cite{ErdI} shows that $m_2(r) \geq 2^{r-1}$ for all $r \in \N^+$. This bound was improved by Beck~\cite{Beck} to $m_2(r) \geq r^{1/3-o(1)} 2^r$, which in particular yields
	\begin{equation}\label{eq:limit}
	\lim_{r \to \infty} \frac{m_2(r)}{2^{r-1}} \,=\, \infty.
	\end{equation}
%	\ep{Spencer presented a simplified proof of Beck's result in \cite{Spencer}.} The best currently known asymptotic lower bound on $m_2(r)$ is due to
	Radhakrishnan and Srinivasan \cite{RS} improved Beck's result further:
	
	\begin{theo}[{Radhakrishnan--Srinivasan \cite{RS}}]\label{theo:RS}
		$m_2(r) = \Omega(\sqrt{r/\log r} \cdot 2^r)$.
	\end{theo}
	\begin{remk}\label{remk:list}
	This is a good place to mention \emph{hypergraph list coloring}, which is a generalization of ordinary hypergraph coloring that is analogous to list coloring of graphs. Since hypergraph list coloring is not the focus of this paper, we shall not define it here; the definition and several relevant results can be found, e.g., in \cite{K_survey}. We point out, however, that \hyperref[theo:RS]{Radhakrishnan and Srinivasan's lower bound} extends, with the same proof, to the list coloring setting \cite[Remark~3]{K_survey}. %Since hypergraph list coloring is not the focus of this paper, we shall not define it here; the definition can be found, e.g., in \cite{K_survey}.
	\end{remk}
	
	The best known asymptotic upper bound on $m_2(r)$ is due to \Erdos \cite{ErdII}:
	
	\begin{theo}[{\Erdos \cite{ErdII}}]\label{theo:Erd_upp}
		$m_2(r) = O(r^2 \cdot 2^r)$.
	\end{theo}

	\Erdos and Lov\'asz \cite{ErdosLovasz} conjectured that neither Theorem~\ref{theo:RS} nor Theorem~\ref{theo:Erd_upp} is sharp, and that the correct order of magnitude is $m_2(r) = \Theta(r \cdot 2^r)$. %The best known asymptotic upper bound on $m_2(r)$ is $m_2(r) \leq O(r^2 2^{r})$ due to \Erdos \cite{ErdII}.
	
	A hypergraph $\H$ is called \emphd{simple}, or \emphd{linear}, if $|e \cap h| \leq 1$ for all distinct $e$, $h \in \H$. The following is the analog of Definition~\ref{defn:ErdHaj} for simple hypergraphs:
	
	\begin{defn}
		For $r$, $k \in \N^+$, let $m_k^\ast(r)$ be the minimum $|\H|$ over all $r$-uniform simple hypergraphs $\H$ without a proper $k$-coloring.
	\end{defn}
	
	\noindent Using the Lov\'asz Local Lemma, \Erdos and Lov\'asz \cite{ErdosLovasz} established the following bounds on $m_k^\ast(r)$:
	
	\begin{theo}[{\Erdos--Lov\'asz \cite{ErdosLovasz}}]\label{theo:EL}
		 For all $k \geq 2$,
		\[
			\lim_{r \to \infty} \frac{\log m_k^\ast(r)}{r} \,=\, 2 \log k.
		\]
		Furthermore, there exist constants $c$, $C > 0$ such that for all $r \in \N^+$,
		\[
			c \cdot r^{-3} \cdot 4^r \,\leq\, m_2^\ast(r) \,\leq\, C \cdot r^4 \cdot 4^r.
		\]
	\end{theo}
	
	The lower bound in Theorem~\ref{theo:EL} was improved by Szab\'o~\cite{Szabo} by a factor of $r^{2-\epsilon}$ and then by Kostochka and Kumbhat~\cite{KK}:

	\begin{theo}[{Kostochka--Kumbhat~\cite{KK}}]\label{theo:KK}
		For each $\epsilon>0$, there is $c > 0$ such that for all $r \geq 2$,
		\[
			m_2^\ast (r) \,\geq\, c \cdot r^{-\epsilon} \cdot 4^r.
		\]
	\end{theo}
	
	 %Notice that the lower bound on $m_2^\ast(r)$ in Theorem~\ref{theo:KK} exceeds the upper bound on $m_2(r)$ given by Theorem~\ref{theo:Erd_upp} by a factor of order $r^{-2-o(1)} 2^r$.
	 
	 \noindent Theorems~\ref{theo:Erd_upp} and \ref{theo:KK} show that $m_2^\ast(r)$ exceeds $m_2(r)$ at least by a factor of order $r^{-2-o(1)} 2^r$.

	\subsection{Hypergraph DP-coloring}
	
	In this paper we introduce and study a generalization of hypergraph coloring, which we call \emph{hypergraph DP-coloring}. 
This is an extension of the similar	
	%It is inspired by the analogous 
	concept for graphs \ep{i.e., $2$-uniform hypergraphs} that was developed by Dvo\v r\'ak and Postle~\cite{DP} under the name \emph{correspondence coloring}. A related notion of \emph{local conflict coloring} was studied independently from the algorithmic point of view by Fraigniaud, Heinrich, and Kosowski~\cite{FHK}.
	
	It will be convenient to identify each function $f$ with its graph, i.e., with the set $\set{(x, y) \,:\, f(x) = y}$, as this allows the use of set-theoretic notation, such as $\subseteq$, $\cap$, $|\cdot|$, etc., for functions. We write $\phi \colon X \rightharpoonup Y$ to indicate that $\phi$ is a partial map from $X$ to $Y$, i.e., a function defined on a subset of $X$ with values in $Y$. The set of all partial maps from $X$ to $Y$ is denoted by $\fun{X}{Y}$.
	
	Let $X$ be a set and let $k \in \N^+$. Consider a family $\F \subseteq \fun{X}{k}$ of partial maps. Since we identify functions with their graphs, each $\phi \in \F$ is a subset of the Cartesian product $X \times k$. Thus, we can view $\F$ as a hypergraph on $X \times k$. In particular, we can apply to $\F$ adjectives such as ``$r$-uniform.'' %; in particular, it makes sense to say that $\F$ is \emphd{$r$-uniform} if $|\phi| = r$ for all $\phi \in \F$.
	The \emphd{domain} $\dom(\F)$ of $\F$ is defined by
	\[
		\dom(\F) \defeq \set{\dom(\phi) \,:\, \phi \in \F}.
	\]
	Thus, $\dom(\F)$ is a hypergraph on $X$. Note that $\F$ is $r$-uniform if and only if so is $\dom(\F)$. %Note that, since we identify functions with their graphs, a 
	
	\begin{defn}
		Let $X$ be a set and let $k \in \N^+$. We say that a function $f \colon X \to k$ \emphd{avoids} a partial map $\phi \colon X \rightharpoonup k$ if $\phi \not \subseteq f$. Given a family $\F \subseteq \fun{X}{k}$ of partial maps, a \emphd{$(k, \F)$-coloring} \ep{or simply an \emphd{$\F$-coloring} if $k$ is understood} is a function $f \colon X \to k$ that avoids all $\phi \in \F$.
	\end{defn}
	%\begin{remk*}
	\noindent Equivalently, viewing $\F$ as a hypergraph on $X \times k$, we can say that an $\F$-coloring is an independent set in $\F$ that intersects every fiber $\set{x} \times k$, $x \in X$, in precisely one vertex.
	%\end{remk*}
	
	\begin{defn}
		Let $\H$ be a hypergraph on a set $X$ and let $k \in \N^+$. A \emphd{$k$-fold cover} of $\H$ \ep{or simply a \emphd{cover} if $k$ is understood} is a family of partial maps $\F \subseteq \fun{X}{k}$ such that:
	\begin{itemize}[label=--]
		\item $\dom(\F) \subseteq \H$; and
		
		\item if $\dom(\phi) = \dom(\psi)$ for some distinct $\phi$, $\psi \in \F$, then $\phi \cap \psi = \0$.
	\end{itemize}
	\end{defn}

	If $\F$ is a $k$-fold cover of a hypergraph $\H$, then for each $e \in \H$, we have
		$|\set{\phi \in \F \,:\, \dom(\phi) = e}| \leq k$.
	In other words, every edge of $\H$ is ``covered'' by at most $k$ pairwise disjoint edges of $\F$.
	
	\begin{defn}\label{defn:DP}
		We say that a hypergraph $\H$ on a set $X$ is \emphd{$k$-DP-colorable} if each $k$-fold cover $\F$ of $\H$ admits a $(k, \F)$-coloring $f \colon X \to k$.
	\end{defn}
	
	A $k$-DP-colorable hypergraph $\H$ is also $k$-colorable in the ordinary sense. To see this, consider the $k$-fold cover of $\H$ formed by the maps of the form $\set{(x, i) \,:\, x \in e}$, where $e \in \H$ and $0 \leq i < k$. The following example illustrates Definition~\ref{defn:DP} and shows that $k$-DP-colorability is a strictly stronger assumption than ordinary $k$-colorability:
	
	\begin{exmp}\label{exmp:K34}
		Let $\mathcal{K}_4^3$ be the complete $3$-uniform hypergraph on $4$ vertices, i.e., the hypergraph on a $4$-element set $V = \set{x_0, x_1, x_2, y}$ given by $\mathcal{K}_4^3 \defeq \set{e \subset V \,:\, |e| = 3}$. Any partition of $V$ into a pair of sets of size $2$ induces a proper $2$-coloring of $\mathcal{K}_4^3$. However, $\mathcal{K}_4^3$ is not $2$-DP-colorable. Indeed, consider the $2$-fold cover $\F = \set{\phi_i,\, \psi_i \,:\, 0 \leq i \leq 3}$ of $\mathcal{K}_3^4$ shown below:
		\begin{table}[H]
				\begin{tabular}{ c || c c c c}
					& $x_0$ & $x_1$ & $x_2$ & $y$ \\\hline\hline
					$\phi_0$ & $0$ & $0$ & $0$ & \\
					$\psi_0$ & $1$ & $1$ & $1$ & \\\hline
					$\phi_1$ & $0$ & $1$ & & $0$ \\
					$\psi_1$ & $1$ & $0$ & & $1$ \\\hline
					$\phi_2$ & $1$ & & $0$ & $0$ \\
					$\psi_2$ & $0$ & & $1$ & $1$ \\\hline
					$\phi_3$ & & $0$ & $1$ & $0$ \\
					$\psi_3$ & & $1$ & $0$ & $1$ 
				\end{tabular}
			\end{table}
		\noindent If $f \colon V \to 2$ is an $\F$-coloring, then, because of $\phi_0$ and $\psi_0$, the colors $f(x_0)$, $f(x_1)$, $f(x_2)$ cannot all be the same. Hence, there exist some $0 \leq i$, $j \leq 2$ such that
		\[
			f(x_i) = 0 \quad \text{and} \quad f(x_{i+1}) = 1, \qquad \text{while} \qquad f(x_j) = 1 \quad \text{and} \quad f(x_{j+1}) = 0.
		\]
		\ep{Addition in the indices is interpreted modulo $3$, i.e., $x_3$ is another name for $x_0$.} If $f(y) = 0$, then \[f \,\supset\, \set{(x_i, 0),\, (x_{i+1}, 1),\, (y, 0)} \,\in\, \set{\phi_1,\, \phi_2,\, \phi_3};\] on  the other hand, if $f(y) = 1$, then \[f \, \supset \,\set{(x_j, 1),\,(x_{j+1}, 0),\, (y, 1)} \in \set{\psi_1,\, \psi_2,\, \psi_3}.\] In both cases, $f$ is not an $\F$-coloring, which is a contradiction. %\hfill $\boxdot$
	\end{exmp}

	\subsection{Results on DP-coloring with 2 colors}
	
	The following are the analogs of $m_2(r)$ and $m_2^*(r)$ in the context of DP-coloring:
	
	\begin{defn}
		For $r$, $k \in \N^+$, let $\tilde{m}_k(r)$ (resp.\ $\tilde{m}_k(r)$) denote the minimum $|\H|$ over all $r$-uniform  (resp.\ $r$-uniform and simple) hypergraphs $\H$ that are not $k$-DP-colorable.
	\end{defn}
	
	\noindent Evidently, $\tilde{m}_k(r) \leq m_k(r)$ and $\tilde{m}^*_k(r) \leq m^*_k(r)$
	 for all $r$, $k \in \N^+$, and it is easy to see that \[\tilde{m}_2(1) = m_2(1) =  \tilde{m}^*_2(1) = m_2^*(1)=1 \qquad \text{and} \qquad \tilde{m}_2(2) = m_2(2) =  \tilde{m}^*_2(2) = m^*_2(2) =3.\] 	
	 
	As with graph DP-coloring, it is interesting to compare the behavior of hypergraph DP-coloring with that of ordinary coloring. Usually, arguments that only rely on the ``local'' structure of the problem transfer easily to the setting of DP-coloring. In particular, the Lov\'asz Local Lemma-based proof of Theorem~\ref{theo:EL} given in \cite{ErdosLovasz} extends verbatim to the DP-coloring framework, so we get the following result ``for free'':

	\begin{theo}\label{theo:ELdp}
		All claims of Theorem~\ref{theo:EL} hold with $\tilde{m}^*_2(r)$ in place of $m_2^\ast(r)$. %\footnote{It might be that the more technical proof of Theorem~\ref{theo:KK} can also be transferred to the DP-coloring setting.}
	\end{theo}

	\noindent However, as we shall see, the situation for general, non-simple, hypergraphs is rather different.
	%On the other hand, 
	
	Example~\ref{exmp:K34} shows that $\tilde{m}_2(3) \leq 4$, while $m_2(3) = 7$. The same straightforward first-moment argument as for ordinary coloring yields
	\begin{equation}\label{eq:silly_lower_bound}
		\tilde{m}_2(r) \geq 2^{r-1} \qquad \text{for all } r \in \N^+,
	\end{equation}
	and hence we in fact have $\tilde{m}_2(3) = 4$. We shall also establish the following bounds:
	\begin{equation}\label{eq:4_and_5}
		9 \leq \tilde{m}_2(4) \leq 10 \qquad\text{and}\qquad 16 \leq \tilde{m}_2(5) \leq 17.
	\end{equation}
	The upper bounds in \eqref{eq:4_and_5} are witnessed by constructions described in Section~\ref{sec:constructions}. The lower bound on $\tilde{m}_2(5)$ is an instance of \eqref{eq:silly_lower_bound}, while the lower bound on $\tilde{m}_2(4)$ follows from a more general result stated below, namely Theorem~\ref{theo:lower_bound}. Potapov \cite{Potapov} has recently improved the upper bound on $\tilde{m}_2(5)$ to $16$ and thus showed that \[\tilde{m}_2(5) \,=\, 16.\] The exact value of $\tilde{m}_2(4)$ remains unknown.
	
	%It is perhaps natural to expect that more sophisticated arguments, like those of Beck, Spencer,  and Radhakrishnan--Srinivasan, could be extended to the realm of DP-coloring to establish a version of \eqref{eq:limit} with $\tilde{m}_2(r)$ in place of $m_2(r)$. However, this is not the case: We show that the ratio $\tilde{m}_2(r)/2^{r-1}$ not only does not approach infinity, but indeed  {converges to $1$}:
	
	Before discussing general bounds on $\tilde{m}_2(r)$, we should mention an important difference between DP-coloring and ordinary coloring. Consider a hypergraph $\H$ on a set $X$ and let $f \colon X \to k$ be a coloring that is \emph{not} proper. Suppose that $e \in \H$ is an $f$-monochromatic edge. We may try to make $e$ non-monochromatic by changing the color of a single vertex $x \in e$, but then some other edge, say $h \in \H$, that was not monochromatic before may become monochromatic. The proof of Theorem~\ref{theo:RS} due to Radhakrishnan and Srinivasan crucially relies on the observation that if this happens, then $e \cap h = \set{x}$, i.e., $e$ and $h$ share only one vertex. The same conclusion would hold even for list coloring; that is why Theorem~\ref{theo:RS} extends to the list coloring setting \ep{see Remark~\ref{remk:list}}. However, the analog of this observation fails in the context of DP-coloring; and indeed, in Section~\ref{sec:upper_bound} we prove the following general upper bound, demonstrating the failure of the \hyperref[theo:RS]{Radhakrishnan--Srinivasan bound} for DP-coloring:
	
	\begin{theo}\label{theo:upper_bound}
		For all $r \in \N^+$, we have
		\[
		\tilde{m}_2(r) \,\leq\, 2^{r-1} + 2^{\left\lfloor r/2 \right\rfloor} \,=\, (1+o(1)) 2^{r-1}.
		\]
	\end{theo}
	
	According to Theorem~\ref{theo:upper_bound}, the asymptotic behavior of $\tilde{m}_2(r)$ is markedly different from that of $m_2(r)$, as, in contrast to \eqref{eq:limit}, the ratio $\tilde{m}_2(r)/2^{r-1}$ converges to $1$ as $r$ tends to infinity. %we have %$\lim_{r \to \infty}\tilde{m}_2(r)/2^{r-1} = 1$.
	%\[\lim_{r \to \infty}\frac{\tilde{m}_2(r)}{2^{r-1}} \,=\, 1.\]
	This makes improving the simple lower bound \eqref{eq:silly_lower_bound} surprisingly challenging. In this direction, we show in Section~\ref{sec:lower_bound} that \eqref{eq:silly_lower_bound} is not sharp for all \emph{even} $r$:
	
	\begin{theo}\label{theo:lower_bound}
		If $r \in \N^+$ is even, then $\tilde{m}_2(r) \geq 2^{r-1} + 1$.
	\end{theo}
	
	By way of contrast, there is little we can say about \emph{odd} $r$. The bound \eqref{eq:silly_lower_bound} \emph{is} trivially sharp for $r = 1$ and, due to Example~\ref{exmp:K34}, for $r = 3$. Potapov \cite{Potapov} also showed that the bound \eqref{eq:silly_lower_bound} is sharp for $r = 5$. Furthermore, he found a way to extend these results to show that
	\[
		\tilde{m}_2(r) \,=\, 2^{r-1} \qquad \text{whenever } r = 3^t5^p \text{ for some } t,\, p \in \N.
	\]
	The following natural questions remain open:
	
	\begin{prob}\label{prob:bounds}
		Do there exist any odd $r \in \N^+$ for which $\tilde{m}_2(r) \geq 2^{r-1} + 1$? Do there exist any $r \in \N^+$ such that $\tilde{m}_2(r) \geq 2^{r-1} + 2$? Is it true that, in fact,
		\[
		\limsup_{r \to \infty} \,(\tilde{m}_2(r) - 2^{r-1}) \,=\, \infty?
		\]
	\end{prob}
	
	\subsection{Further results}\label{subsec:further}
	
	Suppose that $\F \subseteq \fun{X}{2}$ is a finite $r$-uniform family of partial functions that does not admit an $\F$-coloring. If $f \colon X \to 2$ is a uniformly random $2$-coloring, then we have
	\begin{equation}\label{eq:first_moment}
		1 \,=\, \P\left[f \supseteq \phi \text{ for some } \phi \in \F\right] \,\leq\, \sum_{\phi \,\in\, \F} \P\left[f \supseteq \phi\right] \,=\, \frac{|\F|}{2^r}.
	\end{equation}
	Therefore, $|\F| \geq 2^r$. Furthermore, if $|\F| = 2^r$, then the inequality in \eqref{eq:first_moment} is non-strict, which means that the events $\set{f \supseteq \phi}$, $\phi \in \F$, are pairwise incompatible. To put this another way, for any two distinct maps $\phi$, $\psi \in \F$, there is some $x \in \dom(\phi) \cap \dom(\psi)$ such that $\phi(x) \neq \psi (x)$. An obvious example of an $r$-uniform family $\F$ of size $2^r$ with this property is the set of \emph{all} functions from $r$ to $2$. Krotov~\cite{Krotov}, motivated by applications in coding theory, found several constructions of $r$-uniform families $\F \subseteq \fun{(r+1)}{2}$ of size $2^r$ with this property and, in particular, proved that their number is doubly-exponential in $r$.
	
	Say that a family $\F\subseteq \fun{X}{2}$ is \emphd{binary} if for all $e \in \dom(\F)$, we have
	\[|\set{\phi \in \F \,:\, \dom(\phi) = e}| \,\leq\, 2.\]
	In particular, if $\F$ is a $2$-fold cover of a hypergraph on $X$, then $\F$ is binary \ep{but not every binary family is a $2$-fold cover}. In Section~\ref{sec:constructions}, we establish the following fact:
	
	\begin{theo}\label{theo:binary}
		For every $r \in \N^+$, there exists a binary $r$-uniform family $\F \subseteq \fun{\N}{2}$ of size $2^r$ without an $\F$-coloring.
	\end{theo}
	
	Similarly, a family $\F\subseteq \fun{X}{k}$ is \emphd{unary} if for all distinct $\phi$, $\psi \in \F$, we have $\dom(\phi) \neq \dom(\psi)$. Studying $\F$-colorings for unary families $\F$ is related to \emph{graph separation choosability}, introduced by Kratochv\'{\i}l, Tuza, and Voigt in \cite{KTV}. The $2$-uniform case has been recently studied by Dvo\v{r}\'{a}k, Esperet, Kang, and Ozeki \cite{LCC} under the name \emph{least conflict choosability}.
	
	By analogy with $\tilde{m}_k(r)$, we define the following parameter:
	
	\begin{defn}
		For $r \geq 2$, let $\tilde{m}_k'(r)$ denote the smallest cardinality of a unary $r$-uniform family $\F \subseteq \fun{\N}{k}$ without an $\F$-coloring.
	\end{defn}
	
	The value $\tilde{m}'_2(r)$ can be equivalently defined as the smallest number of clauses in an unsatisfiable instance of $r$-\textsf{SAT} where all the clauses use distinct $r$-element sets of variables. We are not aware of any previous study of $\tilde{m}'_2(r)$ from the computer science perspective.
	
	In contrast to Theorem~\ref{theo:binary}, we show that the lower bound $\tilde{m}_2'(r) \geq 2^r$ is never sharp:
	
	\begin{theo}\label{theo:unary_lower}
		For all $r \geq 2$, we have $\tilde{m}_2'(r) \geq 2^r + 1$.
	\end{theo}
	
	\noindent Theorem~\ref{theo:unary_lower} is proved in Section~\ref{sec:lower_bound}, alongside with Theorem~\ref{theo:lower_bound}. We also establish the following upper bound:
	
	\begin{theo}\label{theo:unary_upper}
		For all $r \geq 2$, we have
		$
		\tilde{m}'_2(r) \,\leq\, 2^{r} + 2^{\lceil r/2 \rceil}
		$.
	\end{theo}
	
	\noindent We prove Theorem~\ref{theo:unary_upper} in Section~\ref{sec:upper_bound} and then derive Theorem~\ref{theo:upper_bound} as a corollary.

	\section{Constructions}\label{sec:constructions}
	
	\subsection{4-Uniform hypergraphs}
	
	In this subsection we prove that $\tilde{m}_2(4) \leq 10$; in other words, we construct a non-$2$-DP-colorable $4$-uniform hypergraph with $10$ edges. To begin with, we need a lemma:
	
	\begin{lemma}\label{lemma:neq}
		Let $\mathcal{K}_5^4$ be the complete $4$-uniform hypergraph on a $5$-vertex set $V = \set{x_0, x_1, x_2, y_0, y_1}$.
		
		\begin{enumerate}[label=\ep{\itshape\alph*}, wide]
			\item\label{item:neq} There is a $2$-fold cover $\F_\neq$ of $\mathcal{K}_5^4$ such that if $f \colon V \to 2$ is an $\F_\neq$-coloring, then $f(y_0) \neq f(y_1)$.
			
			\item\label{item:eq} There is a $2$-fold cover $\F_=$ of $\mathcal{K}_5^4$ such that if $f \colon V \to 2$ is an $\F_=$-coloring, then $f(y_0) = f(y_1)$.
		\end{enumerate}
	\end{lemma}
	\begin{coolproof}
		\ref{item:neq} We claim that the cover $\F_\neq = \set{\phi_i,\, \psi_i \,:\, 0 \leq i \leq 4}$ shown in Fig.~\ref{fig:neq} on the left is as desired. 
		\begin{figure}[h]
				\begin{tabular}{ c || c c c c c}
					& $x_0$ & $x_1$ & $x_2$ & $y_0$ & $y_1$ \\\hline\hline
					$\phi_0$ & $0$ & $1$ & & $0$ & $0$ \\
					$\psi_0$ & $1$ & $0$ & & $1$ & $1$ \\\hline
					$\phi_1$ & $1$ & & $0$ & $0$ & $0$ \\
					$\psi_1$ & $0$ & & $1$ & $1$ & $1$ \\\hline
					$\phi_2$ & & $0$ & $1$ & $0$ & $0$ \\
					$\psi_2$ & & $1$ & $0$ & $1$ & $1$ \\\hline
					$\phi_3$ & $0$ & $0$ & $0$ & $0$ & \\
					$\psi_3$ & $1$ & $1$ & $1$ & $1$ & \\\hline
					$\phi_4$ & $0$ & $0$ & $0$ & & $1$ \\
					$\psi_4$ & $1$ & $1$ & $1$ & & $0$ \\
				\end{tabular}
				\qquad\qquad\qquad
				\begin{tabular}{ c || c c c c c}
					& $x_0$ & $x_1$ & $x_2$ & $y_0$ & $y_1$ \\\hline\hline
					$\chi_0$ & $0$ & $1$ & & $0$ & $1$ \\
					$\rho_0$ & $1$ & $0$ & & $1$ & $0$ \\\hline
					$\chi_1$ & $1$ & & $0$ & $0$ & $1$ \\
					$\rho_1$ & $0$ & & $1$ & $1$ & $0$ \\\hline
					$\chi_2$ & & $0$ & $1$ & $0$ & $1$ \\
					$\rho_2$ & & $1$ & $0$ & $1$ & $0$ \\\hline
					$\chi_3$ & $0$ & $0$ & $0$ & $0$ & \\
					$\rho_3$ & $1$ & $1$ & $1$ & $1$ & \\\hline
					$\chi_4$ & $0$ & $0$ & $0$ & & $0$ \\
					$\rho_4$ & $1$ & $1$ & $1$ & & $1$ \\
				\end{tabular}
				\caption{The families $\F_\neq$ \ep{left} and $\F_=$ \ep{right}.}\label{fig:neq}
		\end{figure}
		Indeed, suppose that $f \colon V \to 2$ is an $\F_\neq$-coloring such that $f(y_0) = f(y_1)$. Without loss of generality, assume that $f(y_0) = f(y_1) = 0$. Then, because of $\phi_3$ and $\psi_4$, the colors $f(x_0)$, $f(x_1)$, and $f(x_2)$ cannot all be the same. Hence, there exists some $0 \leq i \leq 2$ with $f(x_i) = 0$ and $f(x_{i+1}) = 1$. \ep{Addition in the indices is interpreted modulo $3$, i.e., $x_3$ is another name for $x_0$.} But then
		\[f \,\supset\, \set{(x_i, 0),\, (x_{i+1}, 1),\, (y_0, 0), \, (y_1, 0)} \,\in\, \set{\phi_0,\, \phi_1,\, \phi_2},\]
		which is a contradiction.
		
		\ref{item:eq} The desired cover is $\F_= = \set{\chi_i,\, \rho_i \,:\, 0 \leq i \leq 4}$ shown in Fig.~\ref{fig:neq} on the right. The proof is analogous to the proof of \ref{item:neq}, and we omit it.
	\end{coolproof}
	
	Now we are ready to prove that $\tilde{m}_2(4) \leq 10$. Let $X_0$, $X_1$, and $Y$ be three pairwise disjoint sets with $|X_0| = |X_1| = 3$ and $|Y| = 2$. For $i \in \set{0,1}$, let $\H_i$ be the complete $4$-uniform hypergraph on $X_i \cup Y$ and let $\H \defeq \H_0 \cup \H_1$. Let the elements of $Y$ be $y_0$ and $y_1$. Using Lemma~\ref{lemma:neq}\ref{item:neq}, we find a $2$-fold cover $\F_0$ of $\H_0$ such that if $f \colon X_0 \cup Y \to 2$ is an $\F_0$-coloring, then $f(y_0)\neq f(y_1)$. Similarly, Lemma~\ref{lemma:neq}\ref{item:eq} gives us a $2$-fold cover $\F_1$ of $\H_1$ such that if $f \colon X_1 \cup Y \to 2$ is an $\F_1$-coloring, then $f(y_0) = f(y_1)$. Then $\F \defeq \F_0 \cup \F_1$ is a $2$-fold cover of $\H$ without an $\F$-coloring. This shows that $\H$ is not $2$-DP-colorable, and hence $\tilde{m}_2(4) \leq |\H| = 10$, as desired.
	
	\subsection{5-Uniform hypergraphs}
	
	In this subsection we prove that $\tilde{m}_2(5) \leq 17$.
	
	\begin{lemma}\label{lemma:nine}
		Let $X = \set{x_0, x_1, x_2}$ and $Y = \set{y_0, y_1, y_2}$ be two disjoint $3$-element sets and let $v$ be an additional vertex not in $X \cup Y$. Let $\H$ be the $5$-uniform hypergraph on $X \cup Y \cup \set{v}$ given by
		\[
			\H \defeq \set{e \subset X \cup Y \cup \set{v} \,:\, |e \cap X| = |e \cap Y| = 2 \text{ and } v \in e}.
		\]
		There exists a $2$-fold cover $\F$ for $\H$ such that if $f \colon X \cup Y \cup \set{v} \to 2$ is an $\F$-coloring, then we have $f(x_0) = f(x_1) = f(x_2)$ or $f(y_0) = f(y_1) = f(y_2)$.
	\end{lemma}
	\begin{coolproof}
		Consider any $e \in \H$. There exist \ep{unique} indices $i$ and $j$ such that
		\[
			e = \set{x_i, \, x_{i+1}, \, y_j, \, y_{j+1}, \, v},
		\]
		where addition is interpreted modulo $3$. Define $\phi_e$, $\psi_e \colon e \to 2$ as follows:
		\[
			\phi_e \defeq \set{(x_i, 0),\, (x_{i+1}, 1), \, (y_j, 0), \, (y_{j+1}, 1), \, (v, 0)},
		\]
		\[
			\psi_e \defeq \set{(x_i, 1),\, (x_{i+1}, 0), \, (y_j, 1), \, (y_{j+1}, 0), \, (v, 1)},
		\]
		and set $\F \defeq \set{\phi_e,\, \psi_e \, :\, e \in \H}$. By construction, $\F$ is a $2$-fold cover of $\H$. Suppose that $f$ is an $\F$-coloring. Without loss of generality, assume that $f(v) = 0$. If neither $f(x_0) = f(x_1) = f(x_2)$ nor $f(y_0) = f(y_1) = f(y_2)$, then we can choose indices $i$ and $j$ so that
		\[
			f(x_i) = f(y_j) = 0 \qquad \text{and} \qquad f(x_{i+1}) = f(y_{j+1}) = 1.
		\]
		But then $\phi_e \subset f$ for $e \defeq \set{x_i, x_{i+1}, y_j, y_{j+1}, v}$; a contradiction. Thus, $\F$ is as desired.
	\end{coolproof}
	
	\begin{lemma}\label{lemma:three}
		Let $X = \set{x_0, x_1, x_2}$ and $Y = \set{y_0, y_1, y_2}$ be two disjoint $3$-element sets and let $\H$ be the $5$-uniform hypergraph on $X \cup Y$ given by
		\[
			\H \defeq \set{e \subset X \cup Y \,:\, |e \cap Y| = 2 \text{ and } X \subset e}.
		\]
		There exists a $2$-fold cover $\F$ of $\H$ such that if $f \colon X \cup Y \to 2$ is an $\F$-coloring with $f(x_0) = f(x_1) = f(x_2)$, then $f(y_0) = f(y_1) = f(y_2)$.
	\end{lemma}
	\begin{coolproof}
		It is straightforward to verify that the family $\F = \set{\phi_i,\, \psi_i \,:\, 0 \leq i \leq 2}$ shown in Fig.~\ref{fig:copy} has all the desired properties.
	\end{coolproof}
	
		\begin{figure}[H]
			\begin{tabular}{ c || c c c c c c}
				& $x_0$ & $x_1$ & $x_2$ & $y_0$ & $y_1$ & $y_2$ \\\hline\hline
				$\phi_0$ & $0$ & $0$ & $0$ & $0$ & $1$ & \\
				$\psi_0$ & $1$ & $1$ & $1$ & $1$ & $0$ & \\\hline
				$\phi_1$ & $0$ & $0$ & $0$ & $1$ & & $0$ \\
				$\psi_1$ & $1$ & $1$ & $1$ & $0$ & & $1$ \\\hline
				$\phi_2$ & $0$ & $0$ & $0$ & & $0$ & $1$ \\
				$\psi_2$ & $1$ & $1$ & $1$ & & $1$ & $0$ \\
			\end{tabular}
			\caption{The family $\F$ in Lemma~\ref{lemma:three}.}\label{fig:copy}
		\end{figure}
	
	Let $X = \set{x_0, x_1, x_2}$, $Y = \set{y_0, y_1, y_2}$, and $Z = \set{z_0, z_1, z_2}$ be three pairwise disjoint $3$-element sets and let $v$ be an additional vertex not in $X \cup Y \cup Z$. We shall build a non-$2$-DP-colorable $5$-uniform hypergraph on $X \cup Y \cup Z \cup \set{v}$ with $17$ edges in four stages. First, let $\H_0$ be the $5$-uniform hypergraph on $X \cup Y \cup \set{v}$ given by
	\[
		\H_0 \defeq \set{e \subset X \cup Y \cup \set{v} \,:\, |e \cap X| = |e \cap Y| = 2 \text{ and } v \in e},
	\]
	i.e., as in Lemma~\ref{lemma:nine}, and let $\F_0$ be a $2$-fold cover of $\H_0$ such that if $f \colon X \cup Y \cup \set{v} \to 2$ is an $\F_0$-coloring, then $f(x_0) = f(x_1) = f(x_2)$ or $f(y_0) = f(y_1) = f(y_2)$. Second, as in Lemma~\ref{lemma:three}, let
	\[
		\H_1 \defeq \set{e \subset X \cup Y \,:\, |e \cap Y| = 2 \text{ and } X \subset e}
	\]
	and let $\F_1$ be a $2$-fold cover of $\H_1$ such that if $f \colon X \cup Y \to 2$ is an $\F_1$-coloring with $f(x_0) = f(x_1) = f(x_2)$, then $f(y_0) = f(y_1) = f(y_2)$. Third, let
	\[
		\H_2 \defeq \set{e \subset Y \cup Z \,:\, |e \cap Z| = 2 \text{ and } Y \subset e}
	\]
	apply Lemma~\ref{lemma:three} with $Y$ in place of $X$ and $Z$ in place of $Y$ to obtain a $2$-fold cover $\F_2$ of $\H_2$ such that if $f \colon Y \cup Z \to 2$ is an $\F_2$-coloring with $f(y_0) = f(y_1) = f(y_2)$, then $f(z_0) = f(z_1) = f(z_2)$. Finally, let
	\[
		\H_3 \defeq \set{\set{y_0, y_1} \cup Z, \, \set{y_1, y_2} \cup Z}
	\]
	and set $\F_3 = \set{\phi_0,\, \psi_0, \, \phi_1, \, \psi_1}$ to be the $2$-fold cover of $\H_3$ shown in Fig.~\ref{fig:two_edges}.
	\begin{figure}[H]
		\begin{tabular}{ c || c c c c c c}
			& $y_0$ & $y_1$ & $y_2$ & $z_0$ & $z_1$ & $z_2$ \\\hline\hline
			$\phi_0$ & $0$ & $0$ & & $0$ & $0$ & $0$ \\
			$\psi_0$ & $1$ & $1$ & & $1$ & $1$ & $1$ \\\hline
			$\phi_1$ & & $1$ & $1$ & $0$ & $0$ & $0$ \\
			$\psi_1$ & & $0$ & $0$ & $1$ & $1$ & $1$ \\
		\end{tabular}
		\caption{The family $\F_3$.}\label{fig:two_edges}
	\end{figure}
	
	Consider the hypergraph $\H \defeq \H_0 \cup \H_1 \cup \H_2 \cup \H_3$. It is $5$-uniform and, since the hypergraphs $\H_0$, $\H_1$, $\H_2$, and $\H_3$ are pairwise edge-disjoint, we have $|\H| = |\H_0| + |\H_1| + |\H_2| + |\H_3| = 9 + 3 + 3 + 2 = 17$. The family $\F \defeq \F_0 \cup \F_1 \cup \F_2 \cup \F_3$ is a $2$-fold cover of $\H$, and we claim that there is no $\F$-coloring. Indeed, suppose that $f \colon X \cup Y \cup Z \cup \set{v} \to 2$ is an $\F$-coloring. Since $f$ is an $\F_0$-coloring, we have $f(x_0) = f(x_1) = f(x_2)$ or $f(y_0) = f(y_1) = f(y_2)$, and since $f$ is also an $\F_1$-coloring, we conclude that $f(y_0) = f(y_1) = f(y_2)$. Since $f$ is an $\F_2$-coloring, this implies $f(z_0) = f(z_1) = f(z_2)$. But then $f$ is not an $\F_3$-coloring---a contradiction.
	
	\subsection{Binary families}\label{subsec:binary}
	
	In this subsection we prove Theorem~\ref{theo:binary}. The proof is by induction on $r \in \N^+$. For $r = 1$, the conclusion of Theorem~\ref{theo:binary} is obvious, so suppose that $r > 1$ and that Theorem~\ref{theo:binary} holds with $r-1$ in place of $r$. For $i \in \set{0,1}$, let $\F_i \subseteq \fun{X_i}{2}$ be a binary $(r-1)$-uniform family of size $2^{r-1}$ without an $\F_i$-coloring, where we assume that $X_0 \cap X_1 = \0$. Let $v$ be an additional vertex not in $X_0 \cup X_1$. Set $X \defeq X_0 \cup X_1 \cup \set{v}$ and define $\F \subseteq \fun{X}{2}$ by
	\[
		\F \defeq \set{\phi \cup \set{(v, 0)}, \, \psi \cup \set{(v, 1)} \,:\, \phi \in \F_0,\, \psi \in \F_1}.
	\]
	Clearly, the family $\F$ is binary and $r$-uniform, and we have $|\F| = 2^r$. It remains to show that there is no $\F$-coloring. To that end, consider any $f \colon X \to 2$. By the choice of $\F_0$ and $\F_1$, there exist some $\phi \in \F_0$ and $\psi \in \F_1$ with $\phi$, $\psi \subset f$. If $f(v) = 0$, then $\phi \cup \set{(v, 0)} \subset f$, while if $f(v) = 1$, then $\psi \cup \set{(v, 1)} \subset f$; in either case, $f$ is not an $\F$-coloring.
	
	\section{Upper bounds}\label{sec:upper_bound}
	
	\subsection{Reduction of Theorem~\ref{theo:upper_bound} to Theorem~\ref{theo:unary_upper}}
	
	%The following piece of notation will be useful: for a function $\phi \colon A \to 2$, let $\overline{\phi} \colon A \to 2$ be given by \[\overline{\phi}(x) \defeq 1 - \phi(x) \quad \text{for all } x \in A.\]
	%Given a set $X$ and a family $\F \subseteq \fun{X}{2}$, let
	%\[
	%	\overline{\F} \defeq \set{\overline{\phi} \,:\, \phi \in \F}.
	%\]
	
	%\begin{lemma}
	%	Let $\F \subseteq \fun{X}{2}$ be a family of partial maps. If there is no $\F$-coloring, then there is also no $\overline{\F}$-coloring.
	%\end{lemma}
	%\begin{coolproof}
	%	If $f \colon X \to 2$ is an $\overline{\F}$-coloring, then $\overline{f}$ is an $\F$-coloring.
	%\end{coolproof}
	
	\begin{lemma}\label{lemma:unary_vs_binary}
		For all $r \geq 3$, we have $\tilde{m}_2(r) \,\leq\, \tilde{m}_2'(r - 1)$ and $\tilde{m}_2'(r) \leq 2 \cdot \tilde{m}_2'(r-1)$.
	\end{lemma}
	\begin{coolproof}
		Set $m \defeq \tilde{m}'_2(r-1)$. Let $\F \subseteq \fun{X}{2}$ be a unary $(r-1)$-uniform family of size $m$ without an $\F$-coloring \ep{such $\F$ exists by the choice of $m$}. Let $v$ be an additional vertex not in $X$ and let
		\[
			\F^\ast \defeq \set{\phi \cup \set{(v, 0)}, \, \overline{\phi} \cup \set{(v, 1)} \,:\, \phi \in \F},
		\]
		where $\overline{\phi} \colon \dom(\phi) \to 2$ is given by $\overline{\phi}(x) \defeq 1 - \phi(x)$.
		Clearly, $\F^\ast$ is a $2$-fold cover of the $r$-uniform hypergraph $\H \defeq \dom(\F^\ast)$ and $|\H| = m$. Suppose that $f \colon X \cup \set{v} \to 2$ is an $\F^\ast$-coloring. Without loss of generality, assume that $f(v) = 0$. By the choice of $\F$, there is some $\phi \in \F$ such that $\phi \subset f$. But then $\phi \cup \set{(v,0)} \subseteq f$, so $f$ is not an $\F^\ast$-coloring---a contradiction. This shows that $\tilde{m}_2(r) \leq m$.
		
		The proof of the inequality $\tilde{m}_2'(r) \leq 2m$ is almost identical to the proof of Theorem~\ref{theo:binary} given in \S\ref{subsec:binary}, so we omit it.
	\end{coolproof}
	
	In view of Lemma~\ref{lemma:unary_vs_binary}, Theorem~\ref{theo:upper_bound} follows from Theorem~\ref{theo:unary_upper} and, furthermore, it is enough to prove Theorem~\ref{theo:unary_upper} for even $r$, which is done in the next subsection.
	
	\subsection{Proof of Theorem~\ref{theo:unary_upper} for even $r$}
	
	\begin{lemma}\label{lemma:construction}
		Let $X = \set{x_0, \ldots, x_{r-1}}$ and $Y = \set{y_0, \ldots, y_{r-1}}$ be two disjoint sets of cardinality $r \in \N^+$. Let $\H$ be the $r$-uniform hypergraph on $X \cup Y$ given by
		\[
			\H \defeq \set{e \subset X \cup Y \,:\, |e \cap \set{x_i, y_i}| = 1 \text{ for all } 0 \leq i \leq r-1}.
		\]
		Then there exists a unary family $\F \subset \fun{X \cup Y}{2}$ such that \[\dom(\F) = \H,\] and every $\F$-coloring $f \colon X \cup Y \to 2$ satisfies the following:
		\begin{itemize}[label=--]
			\item $f(x_i) \neq f(y_i)$ for all $0 \leq i \leq r-1$;
			
			\item the number of indices $i$ such that $f(x_i) = 1$ is odd.
		\end{itemize}
	\end{lemma}
	\begin{coolproof}
		Consider any $e \in \H$. By the definition of $\H$, we can write
		\[
			e = \set{z_0, \ldots, z_{r-1}},
		\]
		where $z_i \in \set{x_i, y_i}$ for each $0 \leq i \leq r-1$. Define a map $\phi_e \colon e \to 2$ as follows:
		\[
			\phi(z_i) \defeq \begin{cases}
				0 &\text{if } z_{i+1} = x_{i+1};\\
				1 &\text{if } z_{i+1} = y_{i+1}. 
			\end{cases}
		\]
		\ep{Addition in the indices is interpreted modulo $r$, i.e., $x_{r}$, $y_{r}$, and $z_{r}$ are other names for $x_0$, $y_0$, and $z_0$ respectively.} Let $\F \defeq \set{\phi_e \,:\, e \in E(\H)}$. Clearly, $\F$ is a unary family such that $\dom(\F) = \H$. We claim that $\F$ is as desired.
		
		Let $f \colon X \cup Y \to 2$ be an arbitrary $\F$-coloring. Construct an edge $e = \set{z_0, \ldots, z_{r-1}} \in \H$ inductively as follows: Set $z_0 \defeq x_0$, and for each $0 \leq i \leq r-2$, let
		\begin{equation}\label{eq:build_edge}
			z_{i+1} \defeq \begin{cases}
				x_{i+1} &\text{if } f(z_{i}) = 0;\\
				y_{i+1} &\text{if } f(z_{i}) = 1.
			\end{cases}
		\end{equation}
		By definition, $\phi_e(z_i) = f(z_i)$ for each $0 \leq i \leq r-2$. But, since $f$ is an $\F$-coloring, we have $\phi_e \not\subset f$. Hence, it must be that $f(z_{r-1}) \neq \phi_e(z_{r-1}) = 0$, i.e., $f(z_{r-1}) = 1$. Therefore, $1 \in \set{f(x_{r-1}), f(y_{r-1})}$. Repeating the same argument starting with $z_0 = y_0$, we also obtain that $0 \in \set{f(x_{r-1}), f(y_{r-1})}$, and thus $f(x_{r-1}) \neq f(y_{r-1})$. Since the family $\F$ is invariant under cyclic permutations of the indices, we conclude that $f(x_i) \neq f(y_i)$ for all $0 \leq i \leq r - 1$.
		
		Now let us have another look at the edge $e = \set{z_0, \ldots, z_{r-1}}$ built using \eqref{eq:build_edge} starting with $z_0 = x_0$. Let $S \defeq \set{0 \leq i \leq r-1 \,:\, f(x_i) = 1}$ and set $S_{i} \defeq S \cap \set{0 \leq j < i}$. We claim that, for each $i$,
		\begin{equation}\label{eq:induction}
			z_i = x_i \quad\Longleftrightarrow\quad |S_{i}| \text{ is even}.
		\end{equation}
		Indeed, \eqref{eq:induction} holds trivially for $i = 0$. Suppose that \eqref{eq:induction} holds for some $i \leq r-2$. Then
		\begin{align*}
			z_{i+1} = x_{i+1} \,&\Longleftrightarrow\, f(z_i) =0 \\
			[\text{since } f(x_i) \neq f(y_i)] \qquad &\Longleftrightarrow\, (z_i = x_i \text{ and } f(x_i) = 0) \text{ or } (z_i = y_i \text{ and } f(x_i) = 1)\\
			&\Longleftrightarrow\, (z_i = x_i \text{ and } i \not \in S) \text{ or } (z_i = y_i \text{ and } i \in S)\\
			[\text{by the inductive hypothesis}] \qquad &\Longleftrightarrow\, (|S_i| \text{ is even and } i \not \in S) \text{ or } (|S_i| \text{ is odd and } i \in S)\\
			&\Longleftrightarrow\, |S_{i+1}| \text{ is even},
		\end{align*}
		as desired. Applying \eqref{eq:induction} with $i = r-1$, we obtain
		\[
			z_{r-1} = x_{r-1} \,\Longleftrightarrow\, |S_{r-1}| \text{ is even}.
		\]
		But we already know that $f(z_{r-1}) = 1$. This means that if $z_{r-1} = x_{r-1}$, then $r-1 \in S$, so $|S| = |S_{r-1}| + 1$ is odd. On the other hand, if $z_{r-1} = y_{r-1}$, then $f(x_{r-1}) = 0$, so $r-1 \not\in S$ and $|S| = |S_{r-1}|$ is odd again. In any case, $|S|$ is odd, as claimed.
	\end{coolproof}
	
	With Lemma~\ref{lemma:construction} in hand, we proceed to prove Theorem~\ref{theo:unary_upper} for even $r$.
	
	\begin{coolproof}[Proof of Theorem~\ref{theo:unary_upper} for even $r$]
		Let $r \in \N^+$ be even and set $\ell \defeq r/2$. Let $X = \set{x_0, \ldots, x_{\ell-1}}$, $Y = \set{y_0, \ldots, y_{\ell-1}}$, $Z = \set{z_0, \ldots, z_{\ell-1}}$, and $W = \set{w_0, \ldots, w_{\ell-1}}$ be four disjoint sets, each of size $\ell$. Let $\H$ be the $r$-uniform hypergraph on $X \cup Y \cup Z \cup W$ given by
		\[
			\H \defeq \set{e \subset X \cup Y \cup Z \cup W \,:\, |e \cap \set{x_i, y_i}| = 1 \text{ and } |e \cap \set{z_i, w_i}| = 1 \text{ for all } 0 \leq i \leq \ell-1}.
		\]
		According to Lemma~\ref{lemma:construction}, there is a unary family $\F \subset \fun{X \cup Y \cup Z \cup W}{2}$ such that $\dom(\F) = \H$, and every $\F$-coloring $f \colon X \cup Y \cup Z \cup W \to 2$ satisfies the following:
		\begin{enumerate}[label={\ep{\emph{\roman*}}}]
			\item\label{item:diff1} $f(x_i) \neq f(y_i)$ and $f(z_i) \neq f(w_i)$ for all $0 \leq i \leq \ell - 1$;
			
			\item\label{item:odd} the cardinality of the set $\set{u \in X \cup Z \,:\, f(u) = 1}$ is odd.
		\end{enumerate}
		Let $\H'$ be the $\ell$-uniform hypergraph on $X \cup W$ given by
		\[
			\H' \defeq \set{e \subset X \cup W \,:\, |e \cap \set{x_i, w_i}| = 1 \text{ for all } 0 \leq i \leq \ell-1},
		\]
		and apply Lemma~\ref{lemma:construction} again to obtain a unary family $\F' \subset \fun{X\cup W}{2}$ such that $\dom(\F') = \H'$, and every $\F'$-coloring $f \colon X \cup W \to 2$ satisfies the following:
		\begin{enumerate}[label={\ep{\emph{\roman*}}},resume]
			\item\label{item:diff2} $f(x_i) \neq f(w_i)$ for all $0 \leq i \leq \ell-1$.%;
			%\item the number of indices $i$ such that $f(x_i) = 1$ is odd.
		\end{enumerate}
		
		\begin{claim*}
			There is no $(\F \cup \F')$-coloring.
		\end{claim*}
		\begin{claimproof}
			Suppose, towards a contradiction, that $f \colon X \cup Y \cup Z \cup W \to 2$ is an $(\F \cup \F')$-coloring. From \ref{item:diff1} and \ref{item:diff2}, it follows that $f(x_i) = f(z_i)$ for all $0 \leq i \leq \ell-1$. But then
			\[
				|\set{u \in X \cup Z \,:\, f(u) = 1}| \,=\, 2 \cdot |\set{x \in X \,:\, f(x) = 1}| \,\equiv\, 0 \pmod 2,
			\]
			contradicting \ref{item:odd}.
		\end{claimproof}
		
		\noindent Note that $|\F \cup \F'| = |\F| + |\F'| = |\H| + |\H'| = 2^r + 2^\ell$, so $\F \cup \F'$ is a unary family of the desired size without a coloring. The only problem is that this family is not $r$-uniform. To fix this, define an $r$-uniform family $\F''$ as follows. For each $\phi \in \F'$, let $\tilde{\phi} \colon X \cup Y \cup Z \cup W \rightharpoonup 2$ be given by
		\[
			\tilde{\phi}(u) \defeq \begin{cases}
				\phi(u) &\text{if } u \in \dom(\phi);\\
				1 - \phi(x_i) &\text{if } x_i \in \dom(\phi) \text{ and } u = y_i;\\
				1 - \phi(w_j) &\text{if } w_j \in \dom(\phi) \text{ and } u = z_j;\\
				\text{undefined} &\text{otherwise}.
			\end{cases}
		\]
		The domain of $\tilde{\phi}$ is obtained from $\dom(\phi)$ by adding all $y_i$ and $z_j$ with $x_i$, $w_j \in \dom(\phi)$. Set \[\F'' \defeq \set{\tilde{\phi} \,:\, \phi \in \F'}.\]
		Then $\F \cup \F''$ is an $r$-uniform unary family of size $2^r + 2^{\ell}$, and we claim that there is no $(\F \cup \F'')$-coloring. Indeed, suppose that $f \colon X \cup Y \cup Z \cup W \to 2$ is an $(\F \cup \F'')$-coloring. Since $f$ cannot be an $\F'$-coloring, there is some $\phi \in \F'$ with $\phi \subset f$. But then, from \ref{item:diff1}, it follows that $\tilde{\phi} \subset f$, which is a contradiction.
	\end{coolproof}
	
	\section{Lower bounds}\label{sec:lower_bound}
	
	\subsection{Families of weight $1$}
	
	We shall deduce Theorems~\ref{theo:lower_bound} and \ref{theo:unary_lower} from a more general statement about families of partial maps that are not necessarily uniform. Before we can state it, we need to introduce a few definitions.
	
	For a map $\phi \colon A \to 2$, where $A$ is a finite set, the \emphd{weight} $\mathsf{w}(\phi)$ of $\phi$ is defined by
	\[
		\mathsf{w}(\phi) \defeq 2^{-|\phi|}.
	\]
	Given a finite set $X$ and a family $\F \subseteq \fun{X}{2}$, define the \emphd{weight} $\mathsf{w}(\F)$ of $\F$ via
	\[
		\mathsf{w}(\F) \defeq \sum_{\phi \in \F} \mathsf{w}(\phi).
	\]
	%We say that $\F$ is \emphd{full} if for every $f \colon X \to \set{0,1}$, there is some $\phi \in \F$ such that $\phi \subseteq f$.
	
	%\noindent The first-moment argument sketched in the beginning of \S\ref{subsec:further} yields the following:
	
	\begin{lemma}\label{lemma:weight_one_prob}
		Let $X$ be a finite set and let $\F \subseteq \fun{X}{2}$ be a family of partial maps without an $\F$-coloring. Then $\mathsf{w}(\F) \geq 1$. Furthermore, if $\mathsf{w}(\F) = 1$, then $|\set{\phi \in \F \,:\, \phi \subseteq f}|=1$ for all $f \in 2^X$.
	\end{lemma}
	\begin{coolproof}
		Pick a coloring $f \in 2^X$ uniformly at random. Then
		\begin{equation}\label{eq:first_moment_1}
		1 \,=\, \P\left[f \supseteq \phi \text{ for some } \phi \in \F\right] \,\leq\, \sum_{\phi \,\in\, \F} \P\left[f \supseteq \phi\right] \,=\, \mathsf{w}(\F),
		\end{equation}
		and hence $\mathsf{w}(\F) \geq 1$. Now suppose that $\mathsf{w}(\F) = 1$. Then the inequality in \eqref{eq:first_moment_1} is non-strict, which implies that  $|\set{\phi \in \F \,:\, \phi \subseteq f}|\leq 1$ for all $f \in 2^X$. On the other hand, we have $|\set{\phi \in \F \,:\, \phi \subseteq f}|\geq 1$ for all $f \in 2^X$ since there is no $\F$-coloring.
	\end{coolproof}
	
	%Our goal in this subsection is to establish a structural property of families $\F \subseteq \fun{X}{2}$ of weight exactly $1$ and with no $\F$-coloring.
	
	\begin{defn}
		Let $X$ be a finite set and let $\F \subseteq \fun{X}{2}$. For $S \subseteq X$, define
	\[
		\F_S \defeq \set{\phi \in \F \,:\, \dom(\phi) \supseteq S},
	\]
	and for $i \in \set{0,1}$, let $\F_S^i$ be the set of all $\phi \in \F_S$ such that $\sum_{x \in S} \phi(x) \equiv i \pmod 2$.
	\end{defn}
	
	%\noindent Notice that the sets $\F_S^0$, $\F_S^1$ form a partition of $\F_S$.
	
	The following lemma is the main result of this subsection:
	
	\begin{lemma}\label{lemma:incexc}
		Let $X$ be a finite set and let $\F \subseteq \fun{X}{2}$. Then, for any $S \subseteq X$, we have
		\begin{equation}\label{eq:parity}
			\mathsf{w}(\F_S^0) - \mathsf{w}(\F_S^1) \,=\, 2^{-|X|} \sum_{f \,\in\, 2^X} (-1)^{\sum_{x \in S} f(x)} |\set{\phi \in \F \,:\, \phi \subseteq f}|.
		\end{equation}
	\end{lemma}
	\begin{coolproof}
		By summing over $\phi \in \F$ first, we can rewrite the right-hand side of \eqref{eq:parity} as follows:
		\[
			2^{-|X|} \sum_{\phi \,\in\, \F} \sum_{f \,\supseteq\, \phi} (-1)^{\sum_{x \in S} f(x)},
		\]
		where the inner summation is over all $f \in 2^X$ with $f \supseteq \phi$.
		
		\begin{claim*}
			Fix any $\phi \colon X \rightharpoonup 2$. Then, for each $S \subseteq X$, we have
			\begin{equation}\label{eq:nonzero} \def\arraystretch{1.5}
				\sum_{f \,\supseteq\, \phi} (-1)^{\sum_{x \in S} f(x)} \,=\, \left\{ \begin{array}{ccl}
				\mathsf{w}(\phi) \cdot 2^{|X|} \cdot (-1)^{\sum_{x \in S} \phi(x)} & \, & \text{if } S \subseteq \dom(\phi);\\
				0 & & \text{if } S \not \subseteq \dom(\phi).
				\end{array}\right.
			\end{equation}
		\end{claim*}
		\begin{claimproof}
			If $S \subseteq \dom(\phi)$, then \[\sum_{x \in S} f(x) = \sum_{x \in S} \phi(x) \qquad \text{for all } f \supseteq \phi,\]
			and, since the number of extensions $f \supseteq \phi$ is equal to $2^{|X| - |\phi|} = \mathsf{w}(\phi) \cdot 2^{|X|}$, we obtain the first case of \eqref{eq:nonzero}. On the other hand, if $S \not \subseteq \dom(\phi)$, then precisely half of the extensions $f \supseteq \phi$ satisfy $\sum_{x \in S} f(x) \equiv 0 \pmod 2$, which yields the second case of \eqref{eq:nonzero}.
		\end{claimproof}
		
		\noindent Using the above claim, we conclude that
		\begin{align*}
			2^{-|X|} \sum_{\phi \,\in\, \F} \sum_{f \,\supseteq\, \phi} (-1)^{\sum_{x \in S} f(x)} \,&=\, \sum_{\phi \,\in\, \F_S} \mathsf{w}(\phi) \cdot (-1)^{\sum_{x\in S} \phi(x)} \\
			%&=\, \sum_{\phi \,\in\, \F_S} \omega(\phi) \cdot (-1)^{\sum_{x\in S} \phi(x)} \\
			&=\, \sum_{\phi \,\in\, \F_S^0} \mathsf{w}(\phi) \,-\, \sum_{\phi \,\in\, \F_S^1} \mathsf{w}(\phi) \,=\, \mathsf{w}(\F_S^0) - \mathsf{w}(\F_S^1), %\qedhere
		\end{align*}
		as desired.
	\end{coolproof}
	
	\begin{corl}\label{corl:weight1}
		Let $X$ be a finite set and let $\F \subseteq \fun{X}{2}$ be a family of weight $1$. If there is no $\F$-coloring, then for every $\0 \neq S \subseteq X$, we have
		\[
			\mathsf{w}(\F_S^0) \,=\, \mathsf{w}(\F_S^1).
		\]
	\end{corl}
	\begin{coolproof}
		Due to Lemma~\ref{lemma:weight_one_prob}, $|\set{\phi \in \F \,:\, \phi \subseteq f}|=1$ for all $f \in 2^X$. Hence, Lemma~\ref{lemma:incexc} yields that for each $S \subseteq X$, we have
		\begin{equation}\label{eq:weight1}
			\mathsf{w}(\F_S^0) - \mathsf{w}(\F_S^1) \,=\, 2^{-|X|} \sum_{f \,\in\, 2^X} (-1)^{\sum_{x \in S} f(x)}.
		\end{equation}
		If $S \neq \0$, then precisely half of the maps $f \in 2^X$ satisfy $\sum_{x \in S} f(x) \equiv 0 \pmod 2$, and hence the right-hand side of \eqref{eq:weight1} is zero, as desired.
	\end{coolproof}
	
	\subsection{Proofs of Theorems~\ref{theo:lower_bound} and \ref{theo:unary_lower}}
	
	For the reader's convenience, we state Theorem~\ref{theo:lower_bound} again:
	
	\begin{theocopy}{theo:lower_bound}
		If $r \in \N^+$ is even, then $\tilde{m}_2(r) \geq 2^{r-1} + 1$.
	\end{theocopy}
	\begin{coolproof}
		Suppose, toward a contradiction, that $\H$ is a non-$2$-DP-colorable $r$-uniform hypergraph with $2^{r-1}$ edges and let $\F$ be a $2$-fold cover of $\H$ without an $\F$-coloring. Then \begin{equation}\label{eq:wub}
		\mathsf{w}(\F) \,=\, |\F| \cdot 2^{-r} \,\leq\, 2 \cdot |\H| \cdot 2^{-r} \,=\, 1,
		\end{equation} and hence in fact $\mathsf{w}(\F) = 1$ by Lemma~\ref{lemma:weight_one_prob}. Therefore, by Corollary~\ref{corl:weight1}, we have $\mathsf{w}(\F_S^0) = \mathsf{w}(\F_S^1)$ for all $S \neq \0$. Consider any $e \in \H$. Since $\mathsf{w}(\F) = 1$, the inequality in \eqref{eq:wub} cannot be strict, and hence $|\F| = 2 |\H|$. This means that $\F_e = \set{\phi, \overline{\phi}}$ for some $\phi \colon e \to 2$ \ep{where $\overline{\phi}$ is the function from $e$ to $2$ given by $\overline{\phi}(x) \defeq 1 - \phi(x)$}. Since $|e| = r$ is even, we have \[\sum_{x \in e} \phi(x) \,\equiv\, \sum_{x \in e} \overline{\phi}(x) \pmod 2,\] and thus exactly one of the sets $\F_e^0$, $\F_e^1$ is nonempty, yielding $\mathsf{w}(\F_e^0) \neq \mathsf{w}(\F_e^1)$; a contradiction.
	\end{coolproof}
	
	A similar argument yields the following strengthening of Theorem~\ref{theo:unary_lower}:
	
	\begin{theo}
		Let $X$ be a finite set and let $\F \subseteq \fun{X}{2}$ be a family of weight $1$ with no $\F$-coloring. Then for every $\phi \in \F$, there is $\psi \in \F \setminus \set{\phi}$ such that $\dom(\psi) \supseteq \dom(\phi)$.
	\end{theo}
	\begin{coolproof}
		By Corollary~\ref{corl:weight1}, we have $\mathsf{w}(\F_S^0) = \mathsf{w}(\F_S^1)$ for all $\0 \neq S \subseteq X$. Suppose that $\phi \in \F$ is such that $\dom(\psi) \not \supseteq \dom(\phi)$ for all $\psi \in \F \setminus \set{\phi}$ and let $S \defeq \dom(\phi)$. Then $\F_S = \set{\phi}$, so exactly one of the sets $\F_S^0$, $\F_S^1$ is nonempty. Hence, $\mathsf{w}(\F_S^0) \neq \mathsf{w}(\F_S^1)$---a contradiction.
	\end{coolproof}
	
	\subsubsection*{{Acknowledgments}}
	
	We are grateful to Michelle Delcourt for pointing out the work of Fraigniaud, Heinrich, and Kosowski \cite{FHK} on local conflict coloring; to Vladimir Potapov for communicating his recent results on DP-coloring of $3^t5^p$-uniform hypergraphs and for drawing our attention to the related research of Krotov \cite{Krotov} in coding theory; and to anonymous referees for their valuable comments.

	{\renewcommand{\markboth}[2]{}% Remove header adjustment
		\printbibliography}
	
\end{document}